\documentclass[twoside,a4paper,10pt]{article}
\usepackage{latexsym,amssymb}
\usepackage{enumerate}
\usepackage{color}
\usepackage{graphicx}
\usepackage{hyperref}
\usepackage[french,english]{babel}
\usepackage[matrix,arrow]{xy}
\def \be{\begin{eqnarray*}}
\def \ee{\end{eqnarray*}}
\def \ben{\begin{enumerate}}
\def \een{\end{enumerate}}
\def \beit{\begin{itemize}}
\def \eeit{\end{itemize}}
\def \bui#1#2{\mathrel{\mathop{\kern 0pt#1}\limits^{#2}}}
\def \buil#1#2{\mathrel{\mathop{\kern 0pt#1}\limits_{#2}}}
\def \bfll{\begin{flushleft}}
\def \efll{\end{flushleft}}
\def \bflr{\begin{flushright}}
\def \eflr{\end{flushright}}

\def \findemo{\hfill$\square$\\}

\def \R{\mathbb{R}}

\def \Ric{\mathrm{Ric}}
\def \ric{\mathrm{ric}}

\newcommand{\rquot}[2]{\raisebox{0.5ex}{$#1$}\!/\!\raisebox{-0.5ex}{$#2$}}

\newtheorem{ethm}{Theorem}[section]

\newtheorem{elemme}[ethm]{Lemma}
\newtheorem{erem}[ethm]{Note}

\newtheorem{eexemple}[ethm]{Example}

\title{A splitting theorem for Riemannian manifolds of generalised Ricci-Hessian type}
\author{Nicolas Ginoux\footnote{Universit\'e de Lorraine, CNRS, IECL, F-57000 Metz, France, E-mail: \texttt{nicolas.ginoux@univ-lorraine.fr}},\, Georges Habib\footnote{Lebanese University, Faculty of Sciences II, Department of Mathematics, P.O. Box 90656 Fanar-Matn, Lebanon,
E-mail: \texttt{ghabib@ul.edu.lb}},\, Ines Kath\footnote{Universit\"at Greifswald, Institut f\"ur Mathematik und Informatik, Walther-Rathenau-Stra\ss{}e 47 17487 Greifswald, Germany, E-mail: \texttt{ines.kath@uni-greifswald.de}
}}

\begin{document}
\maketitle

\noindent\begin{center}\begin{tabular}{p{115mm}}
\begin{small}{\bf Abstract.} In this paper, we study and partially classify those Riemannian manifolds carrying a non-identically vanishing function $f$ whose Hessian is minus $f$ times the Ricci-tensor of the manifold.
\end{small}\\
\end{tabular}\end{center}

$ $\\

\noindent\begin{small}{\it Mathematics Subject Classification} (2010):  53C20, 58J60\\
\end{small}

\noindent\begin{small}{\it Keywords}: Obata equation, Tashiro equation
\end{small}

\section{Introduction}\label{s:intro}

In this paper, we are interested in those Riemannian manifolds $(M^n,g)$ supporting a non-identically-vanishing function $f$ satisfying 
\begin{equation}\label{eq:nabladf=-fric}
\nabla^2f=-f\cdot\Ric
\end{equation}
on $M$, where $\nabla^2f:=\nabla\nabla f$ denotes the Hessian of $f$ and $\Ric$ the Ricci-tensor of $(M^n,g)$, both seen as $(1,1)$-tensor fields.
This equation originates in the search for nontrivial solutions to the so-called \emph{skew-Killing-spinor-equation} \cite{GinouxHabibKathSKS2018}.\\

Equation (\ref{eq:nabladf=-fric}) looks much like that considered by C.~He, P.~Petersen and W.~Wylie {in their search for warped product Einstein metrics} where functions $f$ are considered whose Hessian is a \emph{positive} scalar multiple of $f\cdot(\Ric-\lambda\mathrm{Id})$ for some real constant $\lambda$, see e.g. \cite[Eq. (1.4)]{HePetersenWylie11102456}.
However, no attempt has been made to deal with the \emph{negative} case since.
In another direction, J.~Corvino proved \cite[Prop. 2.7]{Corvino00} that a positive function $f$ satisfies $\nabla^2f=f\cdot\Ric-\left(\Delta f\right)\cdot\mathrm{Id}$ on $M$ if and only if the (Lorentzian) metric $\overline{g}:=-f^2 dt^2\oplus g$ is Einstein on $\R\times M$.
In \cite{HePetersenWylie11102455} the more general situation was considered where the r.h.s. of (\ref{eq:nabladf=-fric}) is replaced by $f\cdot q$ for some {\sl a priori} arbitrary symmetric tensor field $q$ on $TM$; but the statements formulated in \cite{HePetersenWylie11102455} are only valid when, for a given fixed $q$, the space of functions $f$ satisfying $\nabla^2f=f\cdot q$ has dimension at least $2$, see e.g. \cite[Theorem A]{HePetersenWylie11102455}.
Thus we are left with an open problem in case we know about only one such function $f$.\\

{Along the same line, S.~G\"uler and S.A.~Demirba\u{g} define a Riemannian manifold $(M^n,g)$ to be quasi Einstein if and only if there exist smooth functions $u,\alpha,\lambda$ on $M$ such that
$$\Ric+\nabla^2u-\alpha du\otimes\nabla u=\lambda\cdot\mathrm{Id},$$
see \cite[Eq. (1.1)]{GuelerDemirbag18}.
It is easy to see that, for $\alpha=-1$ and $\lambda=0$, a function $u$ solves that equation if and only if $f:=e^u$ solves (\ref{eq:nabladf=-fric}).
However the results of \cite{GuelerDemirbag18} cannot be compared with ours since the quasi Einstein condition seems to be interesting only in the case where $\alpha>0$ and $u>0$.}\\

{ Independently from \cite{HePetersenWylie11102455,HePetersenWylie11102456}, F.E.S.~Feitosa, A.A.~Freitas Filho, J.N.V.~Gomes and R.S.~Pina define gradient almost Ricci soliton warped products by means of functions $f>0,\varphi,\lambda$ satisfying in particular
$$\frac{m}{f}\nabla^2 f+\lambda\cdot\mathrm{Id}=\mathrm{Ric}+\nabla^2\varphi$$
for some nonzero real constant $m$, see \cite[Eq. (1.4)]{FeitosaFreitasFilhoGomesPina15}.
Our equation (\ref{eq:nabladf=-fric}) is the special case of that equation where $m=-1$ and $\lambda=\varphi=0$.
But again \cite{FeitosaFreitasFilhoGomesPina15} only deals with the case where $m>0$; besides, only positive $f$ are considered.
Therefore, no result of \cite{FeitosaFreitasFilhoGomesPina15} can be used in our setting.\\
}

Ricci-flat manifolds carry obvious solutions to (\ref{eq:nabladf=-fric}), just pick constant functions.
Constant functions are actually all solutions to (\ref{eq:nabladf=-fric}) in case $M$ is Ricci-flat and closed.
If $M$ is Ricci-flat, complete but noncompact, then there exists a nonconstant solution to (\ref{eq:nabladf=-fric}) if and only if $M$ is the Riemannian product of $\R$ with a (complete) Ricci-flat manifold $N${; in that case, $f$ is an affine-linear function of the $t\in\R$-coordinate}.
In the search for further nonconstant functions satisfying (\ref{eq:nabladf=-fric}), a natural setting that comes immediately to mind is the case where $(M^n,g)$ is Einstein, because then (\ref{eq:nabladf=-fric}) gets close to the Obata resp. Tashiro equation.
But surprisingly enough only certain two-dimensional spaceforms carry such functions among complete Einstein manifolds, see Lemma \ref{l:eqnabladf=-fric} below.
Further examples can be constructed taking Riemannian products of manifolds carrying a function solving (\ref{eq:nabladf=-fric}) {-- e.g. some $2$-dimensional spaceforms --} with any Ricci-flat Riemannian manifold.
It is however {\sl a priori} unclear whether other examples exist besides those.\\

We show here that, under further geometric assumptions, only products of two-dimensional spaceforms with a Ricci-flat Riemannian manifold can appear, see Theorem \ref{t:eqfRicwitheta}.
In particular, we obtain a classification result covering to some extent the missing case in \cite{HePetersenWylie11102456}.
%The article is organized as follows.

\section{Main result and proof}\label{s:mainresult}

\subsection{Preliminary remarks}\label{ss:prelim}

We start with preliminary results, most of which are elementary or already proved in the literature.
From now on, we shall denote by $S$ the scalar curvature of $M$ and, for any function $h$ on $M$, by $\nabla h$ the gradient vector field of $h$ w.r.t. $g$ on $M$.
First observe that the equation $\nabla^2 f=-f\cdot\Ric$ is of course linear in $f$ but is also invariant under metric rescaling: if $\overline{g}=\lambda^2 g$ for some nonzero real number $\lambda$, then $\overline{\nabla}^2 f=\lambda^{-2}\overline{\nabla}^2 f$ (this comes from the rescaling of the gradient) and $\overline{\Ric}=\lambda^{-2}\Ric$.

\begin{elemme}\label{l:eqnabladf=-fric}
Let $(M^n,g)$ be any connected Riemannian manifold carrying a smooth real-valued function $f$ satisfying {\rm(\ref{eq:nabladf=-fric})} on $M$.
\begin{enumerate}
\item The gradient vector field $\nabla f$ of $f$ w.r.t. $g$ satisfies
\begin{equation}\label{eq:Ricnablaf}
\Ric(\nabla f)=\frac{S}{2}\nabla f+\frac{f}{4}\nabla S.
\end{equation}
\item There exists a real constant $\mu$ such that 
\begin{equation}\label{eq:fLaplacef}
f\Delta f+2|\nabla f|^2=\mu.
\end{equation}
\item If $n>2$ and $f$ is everywhere positive or negative, then $f$ solves {\rm(\ref{eq:nabladf=-fric})} if and only if, setting $u:=\frac{1}{2-n}\ln|f|$, the metric $\overline{g}:=e^{2u}g$ satisfies $\overline{\ric}=(\overline{\Delta}u)\overline{g}-(n-2)(n-3)du\otimes du$ on $M$ and in that case $\overline{\Delta}u=-\frac{\mu}{n-2}e^{2(n-3)u}$.
In particular, if $n=3$, the existence of a positive solution $f$ to {\rm(\ref{eq:nabladf=-fric})} is equivalent to $(M,f^{-2}g)$ being Einstein with scalar curvature $-3\overline{\Delta}\ln|f|$.
\item If $M$ is closed and $f$ is everywhere positive or negative, then $f$ is constant on $M$.
\item If nonempty, the vanishing set $N_0:=f^{-1}(\{0\})$ of $f$ is a scalar-flat totally geodesic hypersurface of $M$.
Moreover, $N_0$ is flat as soon as it is $3$-dimensional and carries a nonzero parallel vector field. 
\item If furthermore $M$ is non-Ricci-flat, Einstein or $2$-dimensional, then $n=2$ and $M$ has constant curvature.
In particular, when $(M^2,g)$ is complete, there exists a nonconstant function $f$ satisfying {\rm(\ref{eq:nabladf=-fric})} if and only if, up to rescaling the metric, the manifold $(M^2,g)$ is isometric to either the round sphere $\mathbb{S}^2$ and $f$ is a nonzero eigenfunction associated to the first positive Laplace eigenvalue; or to flat $\R^2$ or cylinder $\mathbb{S}^1\times\R$ and $f$ is an affine-linear function; or to the hyperbolic plane $\mathbb{H}^2$ and $f$ is a solution to the Tashiro equation $\nabla^2f=f\cdot\mathrm{Id}$.
\item If $S$ is constant, then outside the set of critical points of $f$, the vector field $\nu:=\frac{\nabla f}{|\nabla f|}$ is geodesic.
Moreover, assuming $(M^n,g)$ to be also complete,
\begin{enumerate}
\item if $S>0$, then up to rescaling the metric as well as $f$, we may assume that $S=2$ and that $\mu=f\Delta f+2|\nabla f|^2=2$ on $M$, in which case the function $f$ has $1$ as maximum and $-1$ as minimum value and those are the only critical values of $f$;
\item if $S=0$, then up to rescaling $f$, we may assume that $\mu=2$ on $M$, in which case $f$ has no critical value and $f(M)=\R$, in particular $M$ is noncompact;
\item if $S<0$, then up to rescaling the metric, we may assume that $S=-2$ on $M$, in which case one of the following holds:
\begin{enumerate}
\item if $\mu>0$, then up to rescaling $f$ we may assume that $\mu=2$, in which case $f$ has no critical value and $f(M)=\R$, in particular $M$ is noncompact;
\item if $\mu=0$, then $f$ has no critical value and empty vanishing set and, up to changing $f$ into $-f$, we have $f(M)=(0,\infty)$, in particular $M$ is noncompact;
\item if $\mu<0$, then up to rescaling $f$ we may assume that $\mu=-2$, in which case $f$ has a unique critical value, which, up to changing $f$ into $-f$, can be assumed to be a minimum; moreover, $f(M)=[1,\infty)$, in particular $M$ is noncompact.
\end{enumerate}
\end{enumerate}
\end{enumerate}
\end{elemme}

{\it Proof:}
The proof of the first statement follows that of \cite[Lemma 4]{KimKim03}.
On the one hand, we take the codifferential of $\nabla^2f$ and obtain, choosing a local orthonormal basis $(e_j)_{1\leq j\leq n}$ of $TM$ and using Bochner's formula for $1$-forms:
\begin{eqnarray*}
\delta\nabla^2f&=&-\sum_{j=1}^n\left(\nabla_{e_j}\nabla^2f\right)(e_j)\\
%&=&-\sum_{j=1}^n\nabla_{e_j}\nabla_{e_j}\nabla f-\left(\nabla^2f\right)(\nabla_{e_j}e_j)\\
&=&-\sum_{j=1}^n\nabla_{e_j}\nabla_{e_j}\nabla f-\nabla_{\nabla_{e_j}e_j}\nabla f\\
&=&\nabla^*\nabla (\nabla f)\\
&=&\Delta(\nabla f)-\Ric(\nabla f).
\end{eqnarray*}
On the other hand, by (\ref{eq:nabladf=-fric}) and the formula $\delta\Ric=-\frac{1}{2}\nabla S$,
\begin{eqnarray*}
\delta\nabla^2f&=&\delta\left(-f\cdot\Ric\right)\\
&=&\Ric(\nabla f)-f\cdot\delta\Ric\\
&=&\Ric(\nabla f)+\frac{f}{2}\nabla S.
\end{eqnarray*}
Comparing both identities, we deduce that $\Delta(\nabla f)=2\Ric(\nabla f)+\frac{f}{2}\nabla S$.
But identity (\ref{eq:nabladf=-fric}) also gives 
\begin{equation}\label{eq:DeltaffS}\Delta f=-\mathrm{tr}\left(\nabla^2f\right)=fS,\end{equation}
so that $\Delta(\nabla f)=\nabla (\Delta f)=\nabla(f S)=S\nabla f+f\nabla S$ and therefore $\Ric(\nabla f)=\frac{S}{2}\nabla f+\frac{f}{4}\nabla S$, which is (\ref{eq:Ricnablaf}).\\
By (\ref{eq:nabladf=-fric}) and (\ref{eq:Ricnablaf}), we have
\begin{eqnarray*}
2\nabla(|\nabla f|^2)&=&4\nabla_{\nabla f}^2f\\
&=&-4f\cdot\Ric(\nabla f)\\
&=&-4f\cdot\left(\frac{S}{2}\nabla f+\frac{f}{4}\nabla S\right)\\
&=&-2Sf\nabla f-f^2\nabla S\\
&=&-\nabla(Sf^2)\\
&\bui{=}{\rm(\ref{eq:DeltaffS})}&-\nabla(f\Delta f),
\end{eqnarray*}
% Moreover, (\ref{eq:nabladf=-fric}) implies
% $$-f\cdot\Ric(\nabla f)=\nabla_{\nabla f}^2f=\frac{1}{2}\nabla|\nabla f|^2.$$
% Therefore, the identity $\Delta(\nabla f)=2\Ric(\nabla f)+\frac{f}{2}\nabla S$ yields
% \begin{eqnarray*}
% 0&=&f\cdot\Delta(df)-2f\cdot\Ric(df)-\frac{f^2}{2}d\mathrm{Scal}\\
% &=&f\cdot\left(f\cdot d\mathrm{Scal}+\mathrm{Scal}\cdot df\right)+d|df|^2-\frac{f^2}{2}d\mathrm{Scal}\\
% &=&\frac{f^2}{2}d\mathrm{Scal}+f\cdot\mathrm{Scal}\cdot df+d|df|^2\\
% &=&d\left(\frac{f^2}{2}\cdot\mathrm{Scal}+|df|^2\right)\\
% &=&d\left(\frac{1}{2}f\Delta f+|df|^2\right),
% \end{eqnarray*}
from which (\ref{eq:fLaplacef}) follows.\\
If $f$ vanishes nowhere, then up to changing $f$ into $-f$, we may assume that $f>0$ on $M$.
Writing $f$ as $e^{(2-n)u}$ for some real-valued function $u$ (that is, $u=\frac{1}{2-n}\ln f$), the Ricci-curvatures (as $(0,2
)$-tensor fields) $\ric$ and $\overline{\ric}$ of $(M,g)$ and $(M,\overline{g}=e^{2u}g)$ respectively are related as follows:
\begin{equation}\label{eq:Ricciggbar}
\overline{\ric}=\ric+(2-n)(\nabla du-du\otimes du)+(\Delta u-(n-2)|du|_g^2)g.
\end{equation}
But $\nabla df=(n-2)^2f\cdot du\otimes du+(2-n)f\cdot\nabla du$ and the Laplace operators $\Delta$ of $(M,g)$ and $\overline{\Delta}$ of $(M,\overline{g})$ are related via $\overline{\Delta}v=e^{-2u}\cdot(\Delta v-(n-2)g(du,dv))$ for any function $v$, so that 
\begin{eqnarray*}
\overline{\ric}&=&\ric+\frac{1}{f}\nabla df-(n-2)^2du\otimes du+(n-2) du\otimes du+(\overline{\Delta}u)\overline{g}\\
&=&\ric+\frac{1}{f}\nabla df-(n-2)(n-3)du\otimes du+(\overline{\Delta}u)\overline{g}.
\end{eqnarray*}
As a consequence, $f$ satisfies (\ref{eq:nabladf=-fric}) if and only if $\overline{\ric}=(\overline{\Delta}u)\overline{g}-(n-2)(n-3)du\otimes du$ holds on $M$.
Moreover, 
\begin{eqnarray*}
f\Delta f+2|df|_g^2&=&f\cdot\left(-(n-2)^2f|du|_g^2-(n-2)f\Delta u\right)+2(n-2)^2f^2|du|_g^2\\
&=&-(n-2)f^2\cdot\left(\Delta u-(n-2)|du|_g^2\right)\\
&=&-(n-2)f^2\cdot e^{2u}\cdot\overline{\Delta}u\\
&=&-(n-2)e^{2(2-n)u}\cdot e^{2u}\cdot\overline{\Delta}u\\
&=&-(n-2)e^{2(3-n)u}\cdot\overline{\Delta}u,
%&=&-(n-2)f^{\frac{2(3-n)}{2-n}}\cdot\overline{\Delta}u,
\end{eqnarray*}
in particular (\ref{eq:fLaplacef}) yields $\overline{\Delta}u=-\frac{\mu}{n-2}e^{2(n-3)u}$.
In dimension $3$, we notice that $\overline{\Delta}u=\frac{\overline{S}}{3}$.
This shows the third statement.\\
If $f$ vanishes nowhere, then again we may assume that $f>0$ on $M$.
Since $M$ is closed, $f$ has a minimum and a maximum.
At a point $x$ where $f$ attains its maximum, we have $\mu=f(x)(\Delta f)(x)+2|\nabla_xf|^2=f(x)(\Delta f)(x)\geq0$.
In the same way, $\mu=f(y)(\Delta f)(y)\leq0$ at any point $y$ where $f$ attains its minimum.
We deduce that $\mu=0$ which, by integrating the identity $f\Delta f+2|\nabla f|^2=\mu$ on $M$, yields $df=0$.
This shows the fourth statement.\\
The first part of the fifth statement is the consequence of the following very general fact \cite[Prop. 1.2]{HePetersenWylie11102455}, that we state and reprove here for the sake of completeness: if some smooth real-valued function $f$ satisfies $\nabla^2f=fq$ for some quadratic form $q$ on $M$, then the subset $N_0=f^{-1}\left(\{0\}\right)$ is -- if nonempty -- a totally geodesic smooth hypersurface of $M$.
First, it is a smooth hypersurface because of $d_xf\neq0$ for all $x\in N_0$: namely if $c\colon\mathbb{R}\to M$ is any geodesic with $c(0)=x$, then the function $y:=f\circ c$ satisfies the second order linear ODE $y''=\langle\nabla_{\dot{c}}^2f,\dot{c}\rangle=q(\dot{c},\dot{c})\cdot y$ on $\mathbb{R}$ with the initial condition $y(0)=0$; if $d_xf=0$, then $y'(0)=0$ and hence $y=0$ on $\mathbb{R}$, which would imply that $f=0$ on $M$ by geodesic connectedness, contradiction.
To compute the shape operator $W$ of $N_0$ in $M$, we define $\nu:=\frac{\nabla f}{|\nabla f|}$ to be a unit normal to $N_0$.
Then for all $x\in N_0$ and $X\in T_xM$,
\begin{eqnarray}\label{eq:nablaXnu}
\nonumber\nabla_X^M\mathrm{\nu}&=&X\left(\frac{1}{|\nabla f|}\right)\cdot\nabla f+\frac{1}{|\nabla f|}\cdot\nabla_X^M\nabla f\\
\nonumber&=&-\frac{X\left(|\nabla f|^2\right)}{2|\nabla f|^3}\cdot\nabla f+\frac{1}{|\nabla f|}\cdot\nabla_X^M\nabla f\\
&=&\frac{1}{|\nabla f|}\cdot\left(\nabla_X^2 f-\langle\nabla_X^2 f,\nu\rangle\cdot\nu\right),
\end{eqnarray}
in particular $W_x=-(\nabla\nu)_x=0$ because of $\left(\nabla^2f\right)_x=f(x)q_x=0$.
This shows that $N_0$ lies totally geodesically in $M$.\\
Now recall Gau\ss{} equations for Ricci curvature: for every $X\in TN_0$,
$$\Ric_{N_0}(X)=\Ric(X)^T-R_{X,\nu}^M\nu+\mathrm{tr}_g(W)\cdot WX-W^2X,$$
where $\Ric(X)^T=\Ric(X)-\ric(X,\nu)\nu$ is the component of the Ricci curvature that is tangential to the hypersurface $N_0$.
As a straightforward consequence, if $S_{N_0}$ denotes the scalar curvature of $N_0$,
$$S_{N_0}=S-2\ric(\nu,\nu)+\left(\mathrm{tr}_g(W)\right)^2-|W|^2.$$
Here, $W=0$ and $\Ric(\nu)=\frac{S}{2}\nu$ along $N_0$ because $N_0$ lies totally geodesically in $M$, so that
$$S_{N_0}=S-2\ric(\nu,\nu)=S-S=0.$$
This proves $N_0$ to be scalar-flat.
If moreover $N_0$ is $3$-dimensional and carries a parallel vector field, then it is locally the Riemannian product of a scalar-flat -- and hence flat -- surface with a line, therefore $N_0$ is flat.
This shows the fifth statement.\\
In dimension $2$, we can write $\Ric=\frac{S}{2}\mathrm{Id}=K\mathrm{Id}$, where $K$ is the Gau\ss{} curvature.
But we also know that $\Ric(\nabla f)=\frac{S}{2}\nabla f+\frac{f}{4}\nabla S=K\nabla f+\frac{f}{2}\nabla K$.
Comparing both identities and using the fact that $\{f\neq0\}$ is dense in $M$ leads to $\nabla K=0$, that is, $M$ has constant Gau\ss{} curvature.
Up to rescaling the metric as well as $f$, we may assume that $S,\mu\in\{-2,0,2\}$.
If $M^2$ is complete with constant $S>0$ (hence $K=1$) and $f$ is nonconstant, then $\mu>0$ so that, by Obata's solution to the equation $\nabla^2 f+f\cdot\mathrm{Id}_{TM}=0$, the manifold $M$ must be isometric to the round sphere of radius $1$ and the function $f$ must be a nonzero eigenfunction associated to the first positive eigenvalue of the Laplace operator on $\mathbb{S}^2$, see \cite[Theorem A]{Obata62}.
If $M^2$ is complete and has vanishing curvature, then its universal cover is the flat $\R^2$ and the lift $\tilde{f}$ of $f$ to $\R^2$ must be an affine-linear function of the form $\tilde{f}(x)=\langle a,x\rangle+b$ for some nonzero $a\in\R^2$ and some $b\in\R$; since the only possible nontrivial quotients of $\R^2$ on which $\tilde{f}$ may descend are of the form $\rquot{\R}{\mathbb{Z}\cdot \check{a}}\times\R$ for some nonzero $\check{a}\in a^\perp$, the manifold $M$ itself must be either flat $\R^2$ or such a flat cylinder.
If $M^2$ is complete with constant $S<0$, then $f$ satisfies the Tashiro equation $\nabla^2f=f\cdot\mathrm{Id}_{TM}$.
But then Y.~Tashiro proved that $(M^2,g)$ must be isometric to the hyperbolic plane of constant sectional curvature $-1$, see e.g. \cite[Theorem 2 p.252]{Tashiro65}.
Note that the functions $f$ listed above on $\mathbb{S}^2$, $\mathbb{R}^2$, $\mathbb{S}^1\times\mathbb{R}$ or $\mathbb{H}^2$ obviously satisfy (\ref{eq:nabladf=-fric}).\\
If $(M^n,g)$ is Einstein with $n\geq3$, then it has constant scalar curvature and $\Ric=\frac{S}{n}\cdot\mathrm{Id}$.
%Since $M$ is assumed to be closed, $S$ must be positive (use (\ref{eq:fLaplacef}).
But again the identity $\Ric(\nabla f)=\frac{S}{2}\nabla f+\frac{f}{4}\nabla S=\frac{S}{2}\nabla f$ yields $n=2$ unless $S=0$ and thus $M$ is Ricci-flat.
Therefore, $n=2$ is the only possibility for non-Ricci-flat Einstein $M$.
This shows the sixth statement.\\
If $S$ is constant, then $\Ric(\nabla f)=\frac{S}{2} \nabla f$.
As a consequence, $\nabla_{\nabla f}^2f=-f\Ric(\nabla f)=-\frac{Sf}{2}\nabla f$.
But, as already observed in e.g. \cite[Prop. 1]{RanjSanth95}, away from its vanishing set, the gradient of $f$ is a pointwise eigenvector of its Hessian if and only if the vector field $\nu=\frac{\nabla f}{|\nabla f|}$ is geodesic, see (\ref{eq:nablaXnu}) above.\\
Assuming furthermore $(M^n,g)$ to be complete, we can rescale as before $f$ and $g$ such that $S,\mu\in\{-2,0,2\}$.
In case $S>0$ and hence $S=2$, necessarily $\mu>0$ holds and thus $\mu=2$.
But then $f^2+|\nabla f|^2=1$, so that the only critical points of $f$ are those where $f^2=1$, which by $f^2\leq1$ shows that the only critical points of $f$ are those where $f=\pm1$ and hence where $f$ takes a maximum or minimum value.
Outside critical points of $f$, we may consider the function $y:=f\circ \gamma\colon \R\to\R$, where $\gamma\colon \R\to M$ is a maximal integral curve of the geodesic vector field $\nu$.
Then $y$ satisfies $y'=|\nabla f|\circ\gamma>0$ and $y(t)^2+y'(t)^2=1$, so that $y'=\sqrt{1-y^2}$ and therefore there exists some $\phi\in\R$ such that 
$$y(t)=\cos(t+\phi)\qquad\forall\,t\in \R.$$
Since that function obviously changes sign and $0$ is not a critical value of $f$, we can already deduce that $f$ changes sign, in particular $N_0=f^{-1}(\{0\})$ is nonempty.
Moreover, the explicit formula for $y$ shows that $f$ must have critical points, which are precisely those where $\cos$ reaches its minimum or maximum value.
This shows statement $7.(a)$.\\
In case $S=0$, we have $\Delta f=0$ and therefore (\ref{eq:fLaplacef}) becomes $|\nabla f|^2=\mu$ on $M$, in particular $\mu>0$ and $f$ has no critical point on $M$.
But because of $|\nabla f|=1$, the function $y=f\circ\gamma$ is in fact equal to $t\mapsto t+\phi$ for some constant $\phi\in\R$.
This shows that $f(M)=\R$ and in particular that $M$ cannot be compact.
This proves statement $7.(b)$.\\
In case $S<0$ and thus $S=-2$, there are still three possibilities for $\mu$:
\begin{enumerate}[$\bullet$]
\item If $\mu>0$, then $\mu=2$ and (\ref{eq:fLaplacef}) becomes $-f^2+|\nabla f|^2=1$, hence $f$ has no critical point.
If $\gamma$ is any integral curve of the normalized gradient vector field $\nu=\frac{\nabla f}{|\nabla f|}$, then the function $y:=f\circ\gamma$ satisfies the ODEs $y'=\sqrt{1+y^2}$, therefore $y(t)=\sinh(t+\phi)$ for some real constant $\phi$. 
In particular, $f(M)=\R$ and $M$ cannot be compact.
\item If $\mu=0$, then (\ref{eq:fLaplacef}) becomes $f^2=|\nabla f|^2$.
But since no point where $f$ vanishes can be a critical point by the fifth statement, $f$ has no critical point and therefore must be of constant sign.
Up to turning $f$ into $-f$, we may assume that $f>0$ and thus $f=|\nabla f|$.
Along any integral curve $\gamma$ of $\nu=\frac{\nabla f}{|\nabla f|}$, the function $y:=f\circ\gamma$ satisfies $y'=y$ and hence $y(t)=C\cdot e^t$ for some positive constant $C$.
This shows $f(M)=(0,\infty)$, in particular $M$ cannot be compact.
\item If $\mu<0$, then $\mu=-2$ and (\ref{eq:fLaplacef}) becomes $-f^2+|\nabla f|^2=-1$.
As a consequence, because of $f^2=1+|\nabla f|^2\geq1$, the function $f$ has constant sign and hence we may assume that $f\geq1$ up to changing $f$ into $-f$.
In particular, the only possible critical value of $f$ is $1$, which is an absolute minimum of $f$.
If $\gamma$ is any integral curve of the normalized gradient vector field $\nu=\frac{\nabla f}{|\nabla f|}$, which is defined at least on the set of regular points of $f$, then the function $y:=f\circ\gamma$ satisfies the ODEs $y'=\sqrt{y^2-1}$, therefore $y(t)=\cosh(t+\phi)$ for some real constant $\phi$. 
Since that function has an absolute minimum, it must have a critical point.
It remains to notice that $f(M)=[1,\infty)$ and thus that $M$ cannot be compact.
This shows statement $7.(c)$ and concludes the proof of Lemma \ref{l:eqnabladf=-fric}.
\end{enumerate}

% Namely the covariant derivative of $\nu$ has already been computed above: for all $X\in TM$,
% \begin{eqnarray}\label{eq:nablaXnu}
% \nonumber\nabla_X\nu&=&X(\frac{1}{|\nabla f|})\nabla f+\frac{1}{|\nabla f|}\nabla_X^2f\\
% \nonumber&=&-\frac{X(|\nabla f|^2)}{2|\nabla f|^3}\nabla f+\frac{1}{|\nabla f|}\nabla_X^2f\\
% \nonumber&=&-\frac{g(\nabla_X\nabla f,\nabla f)}{|\nabla f|^3}\nabla f+\frac{1}{|\nabla f|}\nabla_X^2f\\
% \nonumber&=&-\frac{g(\nabla_{X}^2f,\nabla f)}{|\nabla f|^3}\nabla f+\frac{1}{|\nabla f|}\nabla_X^2f\\
% &=&\frac{1}{|\nabla f|}\left(\nabla_X^2f-g(\nabla_X^2f,\nu)\nu\right),
% \end{eqnarray}
% in particular $\nu$ is a pointwise eigenvector of $\nabla^2f$ if and only if $\nabla_\nu\nu=0$.
\findemo

\begin{eexemple}\label{ex:nabladf=-fricdim3}
{\rm In dimension $3$, Lemma \ref{l:eqnabladf=-fric} implies that, starting with any Einstein -- or, equivalently, constant-sectional-curvature- -- manifold $(M^3,g)$ and any real function $u$ such that $\Delta u=\frac{S}{3}$, the function $f:=e^{-u}$ satisfies (\ref{eq:nabladf=-fric}) on the manifold $(M,\overline{g}=e^{-2u}g)$.
In particular, since there is an infinite-dimensional space of harmonic functions on any nonempty open subset $M$ of $\R^3$, there are many nonhomothetic conformal metrics on such $M$ for which nonconstant solutions of (\ref{eq:nabladf=-fric}) exist.
On any nonempty open subset of the $3$-dimensional hyperbolic space $\mathbb{H}^3$ with constant sectional curvature $-1$, there is also an infinite-dimensional affine space of solutions to the Poisson equation $\Delta u=-2${: in geodesic polar coordinates about any fixed point $p\in\mathbb{H}^3$, assuming $u$ to depend only on the geodesic distance $r$ to $p$, that Poisson equation is a second-order linear ODE in $u(r)$ and therefore has infinitely many affinely independent solutions}.
In particular, there are also lots of conformal metrics on $\mathbb{H}^3$ for which nonconstant solutions of (\ref{eq:nabladf=-fric}) exist.\\
{ Note however that, although $\mathbb{H}^3$ is conformally equivalent to the unit open ball $\mathbb{B}^3$ in $\R^3$, we do not obtain the same solutions to the equation depending on the metric we start from.
Namely, we can construct solutions of (\ref{eq:nabladf=-fric}) starting from the Euclidean metric $g$ and from the hyperbolic metric $e^{2w}g$ on $\mathbb{B}^3$, where $e^{2w(x)}=\frac{4}{(1-|x|^2)^2}$ at any $x\in\mathbb{B}^3$.
In both cases we obtain solutions of (\ref{eq:nabladf=-fric}) by conformal change of the metric.
Since $g$ and $e^{2w}g$ lie in the same conformal class, the question arises whether solutions coming from $e^{2w}g$ can coincide with solutions coming from $g$ on $\mathbb{B}^3$.
Assume $f$ were a solution of (\ref{eq:nabladf=-fric}) arising by conformal change of $g$ (by $e^{-2u}$ for some $u\in C^\infty(\mathbb{B}^3)$) and by conformal change of $e^{2w}g$ (by $e^{-2v}$ for some $v\in C^\infty(\mathbb{B}^3)$).
Then $f=e^{-u}=e^{-v}$ and thus $v=u$ would hold, therefore $u$ would satisfy $\Delta_g u=0$ as well as $\Delta_{e^{2w}g}u=-2$, in particular
\begin{eqnarray*}
0&=&\Delta_g u\\
&=&e^{2w}\Delta_{e^{2w}g}u+\langle dw,du\rangle_g\\
&=&-2e^{2w}+\langle dw,du\rangle_g.
\end{eqnarray*}
%where $e^{2w}g$ is the hyperbolic metric.
But the r.h.s. of the last identity has no reason to vanish in general.
Note also that the conformal metrics themselves have no reason to coincide, since otherwise $e^{-2v}e^{2w}g=e^{-2u}g$ would hold hence $u=v-w$ as well and the same kind of argument would lead to an equation that is generally not fulfilled.
}
}
\end{eexemple}

\begin{erem}\label{r:critptf}
\noindent{\rm
% \item In the case where $S$ is constant, the function $f$ is skew-symmetric w.r.t. $N_0$: if $\gamma\colon\,]-\frac{\pi}{2},\frac{\pi}{2}[\to M$ is any integral curve of $\nu$ starting from $N_0$, then $(f\circ\gamma)(-t)=(f\circ\gamma)(t)$ for all $t$ as is obvious from the expression of $f\circ\gamma=\sin$.
{If $S$ is a nonzero constant and $M$ is closed, then the function $f$ is an eigenfunction for the scalar Laplace operator associated to the eigenvalue $S$ on $(M,g)$ and it has at least two nodal domains.
%\footnote{We have not proved that $N_0$ (or any regular hypersurface of $M$) is \emph{connected}.}
Mind however that $S$ is not necessarily the first positive Laplace eigenvalue on $(M,g)$.
E.g. consider the Riemannian manifold $M=\mathbb{S}^2\times\Sigma^{n-2}$ which is the product of standard $\mathbb{S}^2$ with a closed Ricci-flat manifold $\Sigma^{n-2}$, then the first positive Laplace eigenvalue of $\Sigma$ can be made arbitrarily small by rescaling its metric; since the Laplace spectrum of $M$ is the sum of the Laplace spectra of $\mathbb{S}^2$ and $\Sigma$, the first Laplace eigenvalue on $M$ can be made as close to $0$ as desired by rescaling the metric on $\Sigma$.
}
% By the way, does there exist any {\sl a priori} lower eigenvalue bound for the Laplace operator in presence of a nonparallel Killing vector field?
% There is some {\sl \`a la Lichnerowicz} when there is a nontrivial parallel form \cite{Grosjean02}, so why not look for such a bound in our situation?
}
\end{erem}

\subsection{Classification in presence of a particular Killing vector field}\label{ss:withKillingvf}

Next we aim at describing the structure of $M$ using the flow $(F_t^\nu)_t$ of $\nu$.
Namely, outside the possible critical points of $f$, the manifold $M$ is locally diffeomorphic via $(F_t^\nu)_t$ to the product $I\times N_c$ of an open interval with a regular level hypersurface $N_c$ of $f$.
Moreover, the induced metric has the form $dt^2\oplus g_t$ for some one-parameter-family of Riemannian metrics on $N_c$.
To determine $g_t$, one would need to know the Lie derivative of $g$ w.r.t. $\nu$; but for all $X,Y\in TN_c$,
$$(\mathcal{L}_\nu g)(X,Y)=\langle\nabla_X\nu,Y\rangle+\langle\nabla_Y\nu,X\rangle\bui{=}{\rm(\ref{eq:nablaXnu})}\frac{2}{|\nabla f|}\left(\langle\nabla_X^2f,Y\rangle-\langle\nabla_X^2f,\nu\rangle\langle\nu,Y\rangle\right)=-\frac{2f}{|\nabla f|}\ric(X,Y)$$
and we do not know {\sl a priori} more about the Ricci curvature of $M$.
Besides, we have {\sl a priori} no information either on the critical subsets $\{\nabla f=0\}$, we do not even know whether they are totally geodesic submanifolds or not.\\

Therefore, we introduce more assumptions.
We actually introduce some that fit to the particular geometric setting induced by so-called skew Killing spinors, see \cite{GinouxHabibKathSKS2018}.

\begin{ethm}\label{t:eqfRicwitheta}
Let $(M^n,g)$ be a complete Riemannian manifold of dimension $n\geq3$ and constant scalar curvature $S$ carrying a nonconstant real-valued smooth function $f$ satisfying {\rm(\ref{eq:nabladf=-fric})}.
Up to rescaling the metric as well as $f$ we can assume that $S=2\epsilon$ with $\epsilon\in\{-1,0,1\}$ and $\max(f)=1=-\min(f)$ if $S=2$, $|\nabla f|=1$ if $S=0$ and $-f^2+|\nabla f|^2\in\{-1,0,1\}$ if $S=-2$.
Assume also the existence of a non-identically vanishing Killing vector field $\eta$ on $M$ such that
\begin{enumerate}[$\bullet$]
\item the vector fields $\eta$ and $\nabla f$ are pointwise orthogonal,
\item the vector field $\eta_{|_{N_c}}$ is parallel along $N_c:=f^{-1}(\{c\})$ for every regular value $c$ of $f$,
\item if $S\neq0$, the vector field $\eta$ satisfies $\buil{\inf}{M}(|\eta|)=\buil{\inf}{M}(|\nabla f|)$ {and vanishes where $\nabla f$ does},
\item also $\nabla_\eta\eta=\epsilon f\nabla f$ holds in case $n>4$ or $f$ has no critical point.
\end{enumerate}
Then { $M$ is isometric to either the Riemannian product $S(\epsilon)\times\Sigma^{n-2}$ of the simply-connected complete surface with curvature $\epsilon$ with a complete Ricci-flat manifold $\Sigma$ in case $\epsilon\neq0$, or to the Riemannian product of $\R$ with some complete Ricci-flat manifold carrying a nonzero parallel vector field in case $\epsilon=0$.}
%and $f$ comes from an eigenfunction associated to the first positive eigenvalue of the Laplace operator on $\mathbb{S}^2$.
\end{ethm}

{\it Proof:}
The assumption that $\eta\perp\nabla f$ not only means that the flow of $\eta$ preserves the level hypersurfaces of $f$, but also implies that $[\eta,\nabla f]=[\eta,\nu]=0$: for 
$$0=X(\langle\eta,\nabla f\rangle)=\langle\nabla_X\eta,\nabla f\rangle+\langle\eta,\nabla_X^2f\rangle=\langle\nabla_\eta^2f-\nabla_{\nabla f}\eta,X\rangle\qquad\forall\,X\in TM,$$
so that $\nabla_\eta^2f=\nabla_{\nabla f}\eta$, that is, $\nabla_{\nabla f}\eta=-f\Ric(\eta)$ and it also follows that
$$[\eta,\nabla f]=\nabla_\eta\nabla f-\nabla_{\nabla f}\eta=\nabla_\eta^2f-\nabla_\eta^2f=0.$$
As a further consequence of $[\eta,\nabla f]=0$, using again that $\eta$ is Killing,
$$[\eta,\nu]=\eta(\frac{1}{|\nabla f|})\nabla f=-\frac{\langle\nabla_\eta\nabla f,\nabla f\rangle}{|\nabla f|^3}\nabla f=-\frac{\langle\nabla_{\nabla f}\eta,\nabla f\rangle}{|\nabla f|^3}\nabla f=0.$$
In particular, the flow of $\nu$ preserves $\eta$ and conversely the flow of $\eta$ preserves both $\nabla f$ and $\nu$.\\
Next we examine the assumption that $\eta_{|_{N_c}}$ is parallel on $N_c=f^{-1}(\{c\})$, which is a smooth hypersurface for all but finitely many values of $c$ by Lemma \ref{l:eqnabladf=-fric}.
By Gau\ss{}-Weingarten formula, $\nabla^{N_c}\eta=0$ is equivalent to 
$$\nabla _X\eta=\nabla_X^{N_c}\eta+\langle W\eta,X\rangle\nu=\langle W\eta,X\rangle\nu\bui{=}{\rm(\ref{eq:nablaXnu})}-\frac{1}{|\nabla f|}\langle\nabla_\eta^2f,X\rangle\nu=\frac{f}{|\nabla f|}\langle\Ric(\eta),X\rangle\nu$$
for all $X\in TN_c$, where $W=-\nabla \nu$ is the Weingarten-endomorphism-field of $N_c$ in $M$.
With $\nabla_\nu\eta=-\frac{f}{|\nabla f|}\Ric(\eta)$, the above identity is equivalent to
\begin{equation}\label{eq:nablaeta}\nabla \eta=\frac{f}{|\nabla f|}\cdot\left(\Ric(\eta)\otimes\nu-\nu\otimes\Ric(\eta)\right)=\frac{f}{|\nabla f|}\cdot\Ric(\eta)\wedge\nu.\end{equation}
In particular $\mathrm{rk}(\nabla \eta)\leq2$ on the subset of regular points of $f$.
Moreover, $\eta$ cannot vanish anywhere on the subset of regular points of $f$: for if $\eta$ vanished at some regular point $x$, then $\eta$ would vanish along the level hypersurface containing $x$ and, being preserved by the flow of $\nu$, it would have to vanish identically on a nonempty open subset of $M$ and therefore on $M$, which would be a contradiction.
Thus $\eta^{-1}(\{0\})\subset(\nabla f)^{-1}(\{0\})$.
On the other hand, the assumption { $(\nabla f)^{-1}(\{0\})\subset\eta^{-1}(\{0\})$ yields $\eta^{-1}(\{0\})=(\nabla f)^{-1}(\{0\})$}.
%\footnote{The conclusion holds iff $\eta$ vanishes where $\nabla f$ vanishes.}
In particular, when nonempty, $(\nabla f)^{-1}(\{0\})$ is a totally geodesic submanifold of $M$ (vanishing set of a Killing vector field) and has even codimension, which is positive otherwise $\eta$ would vanish identically.\\
By e.g. \cite[Sec. 2.5]{Kocomplexgeom}, the tangent bundle of $\eta^{-1}(\{0\})$ is given by $\ker(\nabla \eta)$ and therefore it has pointwise dimension at most $n-2$.
But since $\mathrm{rk}(\nabla \eta)\leq2$ on $M\setminus\eta^{-1}(\{0\})$, which is a dense open subset of $M$, the inequality $\mathrm{rk}(\nabla \eta)\leq2$ must hold along $\eta^{-1}(\{0\})$ by continuity, in particular $\dim(\ker(\nabla \eta))\geq n-2$ and thus $\dim(\ker(\nabla \eta))=n-2$ along $\eta^{-1}(\{0\})$.
On the whole, when nonempty, the set of critical points of $f$ is a possibly disconnected $(n-2)$-dimensional totally geodesic submanifold of $M$.\\
As a further step, we translate Gau\ss{} equations for Ricci curvature along each $N_c$ in our context.
Denoting $W=-\nabla \nu=\frac{f}{|\nabla f|}\Ric^T$ the Weingarten-endomorphism-field of $N_c$ in $M$, where $\Ric^T$ is the pointwise othogonal projection of $\Ric$ onto $TM$, we have $\mathrm{tr}(W)=\frac{f}{|\nabla f|}\cdot\frac{S}{2}$ by $\Ric(\nu)=\frac{S}{2}\nu$.
As a consequence, we have, for all $X\in TN_c$:
\begin{eqnarray*}
\Ric(X)&=&\Ric(X)^T\\
&=&\Ric_{N_c}(X)+W^2X-\mathrm{tr}(W)WX+R_{X,\nu}\nu\\
&=&\Ric_{N_c}(X)+\frac{f^2}{|\nabla f|^2}\cdot\left(\Ric^2(X)-\frac{S}{2}\Ric(X)\right)+R_{X,\nu}\nu.
\end{eqnarray*}
We can compute the curvature term $R_{X,\nu}\nu$ explicitely: choosing a vector field $X$ that is pointwise tangent to the regular level hypersurfaces of $f$, we can always assume w.l.o.g. that $[X,\nu]=0$ at the point where we compute, so that, with $\nabla_\nu\nu=0$,
\begin{eqnarray*}
R_{X,\nu}\nu&=&-\nabla_\nu\nabla_X\nu\\
&=&\nabla_\nu\left(\frac{f}{|\nabla f|}\Ric(X)\right)\\
&=&\nu\left(\frac{f}{|\nabla f|}\right)\cdot\Ric(X)+\frac{f}{|\nabla f|}\cdot(\nabla_\nu\Ric)(X)+\frac{f}{|\nabla f|}\cdot\Ric(\nabla_\nu X)\\
&=&\frac{\nu(f)|\nabla f|-f\nu(|\nabla f|)}{|\nabla f|^2}\cdot\Ric(X)+\frac{f}{|\nabla f|}\cdot(\nabla_\nu\Ric)(X)+\frac{f}{|\nabla f|}\cdot\Ric(\nabla_X\nu)\\
&=&\frac{|\nabla f|^2-f\frac{\langle\nabla_\nu\nabla f,\nabla f\rangle}{|\nabla f|}}{|\nabla f|^2}\cdot\Ric(X)+\frac{f}{|\nabla f|}\cdot(\nabla_\nu\Ric)(X)-\frac{f^2}{|\nabla f|^2}\cdot\Ric^2(X)\\
&=&\left(1+\frac{Sf^2}{2|\nabla f|^2}\right)\cdot\Ric(X)+\frac{f}{|\nabla f|}\cdot(\nabla_\nu\Ric)(X)-\frac{f^2}{|\nabla f|^2}\cdot\Ric^2(X).
\end{eqnarray*}
We can deduce that, for every $X\in TN_c$,
\begin{eqnarray}\label{eq:RicNc}
\nonumber\Ric_{N_c}(X)&=&\Ric(X)-\frac{f^2}{|\nabla f|^2}\cdot\left(\Ric^2(X)-\frac{S}{2}\Ric(X)\right)-R_{X,\nu}\nu\\
\nonumber&=&\Ric(X)-\frac{f^2}{|\nabla f|^2}\cdot\left(\Ric^2(X)-\frac{S}{2}\Ric(X)\right)\\
\nonumber&&-\left(1+\frac{Sf^2}{2|\nabla f|^2}\right)\cdot\Ric(X)-\frac{f}{|\nabla f|}\cdot(\nabla_\nu\Ric)(X)+\frac{f^2}{|\nabla f|^2}\cdot\Ric^2(X)\\
&=&-\frac{f}{|\nabla f|}\cdot(\nabla_\nu\Ric)(X).
\end{eqnarray}
That identity has important consequences.
First, choosing a local o.n.b. $(e_j)_{1\leq j\leq n-1}$ of $TN_c$,
\begin{eqnarray*}
S_{N_c}&=&\sum_{j=1}^n\langle\Ric_{N_c}(e_j),e_j\rangle\\
&=&-\frac{f}{|\nabla f|}\cdot\sum_{j=1}^n\langle(\nabla_\nu\Ric)(e_j),e_j\rangle\\
&=&-\frac{f}{|\nabla f|}\cdot\left(\sum_{j=1}^n\langle(\nabla_\nu\Ric)(e_j),e_j\rangle+\langle(\nabla_\nu\Ric)(\nu),\nu\rangle\right)+\frac{f}{|\nabla f|}\cdot \langle(\nabla_\nu\Ric)(\nu),\nu\rangle,
\end{eqnarray*}
with $(\nabla_\nu\Ric)(\nu)=\nabla_\nu(\Ric(\nu))-\Ric(\nabla_\nu\nu)=\nabla_\nu(\frac{S}{2}\nu)=0$, so that 
$$S_{N_c}=-\frac{f}{|\nabla f|}\cdot\mathrm{tr}(\nabla_\nu\Ric)=-\frac{f}{|\nabla f|}\cdot\nu(\mathrm{tr}(\Ric))=-\frac{f}{|\nabla f|}\cdot\nu(S)=0.$$
Therefore, each level hypersurface $N_c$ is scalar-flat.\\
In the case where $n=3$ or $4$, the manifold $N_c$ is locally the Riemannian product of a flat manifold with an interval and is hence also flat, in particular $\Ric_{N_c}=0$, which in turn implies that
\begin{equation}\label{eq:nablanuRic}
\nabla_\nu\Ric=0.
\end{equation}
This equation, which holds on the dense open subset $\{\nabla f\neq0\}$, means that all eigenspaces and eigenvalues of the Ricci-tensor of $M$ are preserved under parallel transport along integral curves of $\nu$.
Assume first that $f$ has critical points, in particular $S\neq0$.
Along the critical submanifold $N_{\rm crit}:=(\nabla f)^{-1}(\{0\})$, one has $\ker(\Ric)\supset TN_{\rm crit}$: if $c\colon I\to  N_{\rm crit}$ is any smooth curve, then $f\circ c$ is constant and therefore $0=(f\circ c)''=\langle\nabla_{\dot{c}}^2f,\dot{c}\rangle$ (the gradient of $f$ vanishes along $N_{\rm crit}$), so that $\ric(\dot{c},\dot{c})=0$.
But $\nabla^2 f$ and thus $\Ric$ is either nonpositive or nonnegative along $N_{\rm crit}$ because $N_{\rm crit}$ is a set of minima or maxima of $f$ as we have seen in Lemma \ref{l:eqnabladf=-fric}, therefore $\Ric(\dot{c})=0$.
In particular, $0$ is an eigenvalue of multiplicity at least $n-2$ of the Ricci-tensor; since the Ricci-eigenvalues are constant along the integral curves of $\nu$, it can be deduced that $0$ is an eigenvalue of multiplicity at least $n-2$ everywhere in $M$.
But the multiplicity cannot be greater that $n-2$, otherwise $\Ric$ would have only one nonzero eigenvalue (namely $\frac{S}{2}\in\{\pm1\}$) and hence its trace would be $\frac{S}{2}$, contradiction.
Therefore $0$ is an eigenvalue of multiplicity exactly $n-2$ of $\Ric$ at every point in $M$.
It remains to notice that at regular points, one has $\Ric(\nu)=\frac{S}{2}\nu$ and $\Ric(\eta)\perp\nu$, so that, using $\Ric(\eta)\perp\ker(\Ric)$, we deduce that $\Ric^2(\eta)$ is proportional to $\Ric(\eta)$, the eigenvalue being necessarily equal to $\frac{S}{2}$, that is,
\begin{equation}\label{eq:Ric2eta}
\Ric^2(\eta)=\frac{S}{2}\Ric(\eta).
\end{equation}
This allows for $\eta$ to be normalized as we explain next.
Namely we would like $\nabla_\eta\eta=\epsilon f\nabla f$ to hold on $M$.
Let $\gamma\colon\,(-\varepsilon,\varepsilon)\to M$ be any integral curve of $\nu$ with starting point $\gamma(0)$ in some regular level hypersurface of $f$; we have already seen in the proof of Lemma \ref{l:eqnabladf=-fric} that $y:=f\circ\gamma$ does not depend on the starting point $\gamma(0)$ of $\gamma$ in a fixed level hypersurface of $f$.
Since, as explained above, the vector field $\nu$ is geodesic on $M$ and $\eta$ is parallel along each $N_c$, the function $h:=|\eta|^2\circ\gamma$ only depends on $t$ and not on the starting point $\gamma(0)$.
In other words, $\nabla(|\eta|^2)=\nu(|\eta|^2)\nu$.
But 
\begin{equation}\label{eq:nuetaeta}\nu(|\eta|^2)=2\langle\nabla_\nu\eta,\eta\rangle=-2\langle\nabla_\eta\eta,\nu\rangle=-2\frac{f}{|\nabla f|}\ric(\eta,\eta),\end{equation}
so that $\ric(\eta,\eta)$ also only depends on $t$.
By (\ref{eq:RicNc}) and $\Ric_{N_c}(\eta)=0$ because of $\eta_{|_{N_c}}$ being parallel, we have, outside $N_0=f^{-1}(\{0\})$,
\begin{eqnarray}\label{eq:nuricetaeta}
\nonumber\nu(\ric(\eta,\eta))&=&(\nabla_\nu\ric)(\eta,\eta)+2\ric(\nabla_\nu\eta,\eta)\\
\nonumber&\bui{=}{\rm(\ref{eq:RicNc})}&-\frac{|\nabla f|}{f}\underbrace{\langle\Ric_{N_c}(\eta),\eta\rangle}_{0}-2\frac{f}{|\nabla f|}\ric(\Ric(\eta),\eta)\\
&=&-2\frac{f}{|\nabla f|}|\Ric(\eta)|^2.
\end{eqnarray}
Note here that both (\ref{eq:nuetaeta}) and (\ref{eq:nuricetaeta}) are valid in any dimension and without the condition (\ref{eq:Ric2eta}).
Combining (\ref{eq:nuricetaeta}) with (\ref{eq:Ric2eta}) and using $\nu(|\nabla f|^2)=-2\epsilon f|\nabla f|$, we deduce that 
$$\nu(\ric(\eta,\eta))=-2\epsilon\frac{f}{|\nabla f|}\ric(\eta,\eta)=\frac{\nu(|\nabla f|^2)}{|\nabla f|^2}\ric(\eta,\eta).$$
As a consequence, there exists a real constant $C$, that has the sign of $\epsilon$, such that $\ric(\eta,\eta)=C\cdot|\nabla f|^2$ and thus $\nabla_\eta\eta=Cf\nabla f$ on $M$.
Therefore, up to replacing $\eta$ by $\frac{1}{\sqrt{C\epsilon}}\cdot\eta$, we may assume that $\nabla_\eta\eta=\epsilon f\nabla f$ on $M$.
Note that this concerns only the case where $n\in\{3,4\}$ and $f$ has critical points, otherwise we \emph{assume} $\nabla_\eta\eta=\epsilon f\nabla f$ to hold on $M$.\\
% As a consequence, the first derivative of $h$ also only depends on $t$.
% But 
% $$h'=2\langle\nabla_{\dot{\gamma}}\eta,\eta\rangle=-2(\frac{f}{|\nabla f|}\ric(\eta,\eta))\circ\gamma=-2\frac{\sin}{\cos}\cdot(\ric(\eta,\eta)\circ\gamma),$$
% in particular $\ric(\eta,\eta)\circ\gamma$ also only depends on $t$ and not on $\gamma(0)$.
% Now 
% \begin{eqnarray*}
% \left(\ric(\eta,\eta)\circ\gamma\right)'&=&(\nabla_{\dot{\gamma}}\ric)(\eta,\eta)+2\ric(\nabla_{\dot{\gamma}}\eta,\eta)\\
% &\bui{=}{\rm(\ref{eq:nablanuRic})}&\underbrace{(\nabla_{\nu}\ric)(\eta,\eta)}_{0}+2\ric(\nabla_{\nu}\eta,\eta)\\
% &=&-2\frac{f}{|\nabla f|}\circ\gamma\cdot\ric(\Ric(\eta),\eta)\circ\gamma\\
% &\bui{=}{\rm(\ref{eq:Ric2eta})}&-2\frac{\sin}{\cos}\cdot\frac{S}{2}\cdot(\ric(\eta,\eta)\circ\gamma)\\
% &=&-2\frac{\sin}{\cos}\cdot(\ric(\eta,\eta)\circ\gamma),
% \end{eqnarray*}
% so that $\ric(\eta,\eta)\circ\gamma=C\cdot\cos^2$ for some necessarily positive constant $C$.
% It follows that $\ric(\eta,\eta)=C\cdot|\nabla f|^2$ and hence $\nabla_{\eta}\eta=\frac{f}{|\nabla f|}\ric(\eta,\eta)\nu=Cf\nabla f$ on $M$.
% Therefore, up to replacing $\eta$ by $\frac{1}{\sqrt{C}}\eta$, we can assume that $\nabla_\eta\eta=f\nabla f$ holds on $M$.
Assuming from now on $\nabla_\eta\eta=\epsilon f\nabla f$ and $n\geq3$ to hold, it can be deduced that $|\eta|^2=-\epsilon f^2+\mathrm{cst}$ for some $\mathrm{cst}\in\R$:
Namely $\nabla(|\eta|^2)=-2\nabla_\eta\eta=-2\epsilon f\nabla f=-\epsilon\nabla(f^2)$ and the set of regular points of $f$ is connected.
Moreover, using e.g. (\ref{eq:nuetaeta}), we have $\ric(\eta,\eta)=\epsilon|\nabla f|^2$; differentiating that identity w.r.t. $\nu$ and using (\ref{eq:nuricetaeta}) yields $|\Ric(\eta)|=|\epsilon|\cdot|\nabla f|$.
In case $S=2$, we have $\eta=0$ on $(\nabla f)^{-1}(\{0\})=f^{-1}(\{\pm1\})$, so that $\mathrm{cst}=1$ and hence $|\eta|=\sqrt{1-f^2}=|\nabla f|$.
In case $S=0$, we have $|\eta|^2=\mathrm{cst}$, from which $\ric(\eta,\eta)=0$ and even $\Ric(\eta)=0$ follow.
{In case $S=-2$, we have $-f^2+|\nabla f|^2=\frac{\mu}{2}\in\{-1,0,1\}$, so that $|\eta|^2=f^2+\mathrm{cst}=|\nabla f|^2-\frac{\mu}{2}+\mathrm{cst}$.
Now the assumption $\buil{\inf}{M}(|\eta|)=\buil{\inf}{M}(|\nabla f|)$ yields $\mathrm{cst}-\frac{\mu}{2}=0$, in particular $|\eta|=|\nabla f|$.\\}
% we look at the three possible cases according to the value of $\mu$.
% If $\mu=2$, then $|\eta|^2=f^2+\mathrm{cst}=|\nabla f|^2-1+\mathrm{cst}$.
% The assumption $\buil{\inf}{M}(|\eta|)=\buil{\inf}{M}(|\nabla f|)$ together with $\inf_M(|\nabla f|)=1$ ensures that $\mathrm{cst}=1$, so that $|\eta|=|\nabla f|$.
% If $\mu=0$, then $|\eta|^2=\mathrm{cst}+|\nabla f|^2$ and again the assumption $\buil{\inf}{M}(|\eta|)=\buil{\inf}{M}(|\nabla f|)$ implies $\mathrm{cst}=0$, so that $|\eta|=|\nabla f|$.
% If $\mu=-2$, then $|\eta|^2=1+|\nabla f|^2+\mathrm{cst}$; with $\eta^{-1}(\{0\})=(\nabla f)^{-1}(\{0\})$, we obtain $\mathrm{cst}=-1$ and hence $|\eta|=|\nabla f|$.\\
To sum up, in all cases we obtain $\ric(\eta,\eta)=\epsilon|\Ric(\eta)|\cdot|\eta|$, which is exactly the equality case in Cauchy-Schwarz inequality.
We can thus deduce that $\Ric(\eta)$ is proportional to $\eta$ and hence $\Ric(\eta)=\frac{S}{2}\eta$.\\

If $n>4$ or if $f$ has no critical point, it remains to show that $\Ric(\eta)=\frac{S}{2}\eta$ implies $\ker(\Ric)=\{\eta,\nu\}^\perp$: for we already know from $\Ric_{N_c}=-\frac{f}{|\nabla f|}\cdot(\nabla_\nu\Ric)$ that $S_{N_c}=0$.
But by the Gau\ss{} formula, $S_{N_c}=S-2\ric(\nu,\nu)+\mathrm{tr}(W)^2-|W|^2$, so that, with $S-2\ric(\nu,\nu)=0$, we deduce that 
\begin{eqnarray*}
0&=&\frac{f^2}{|\nabla f|^2}\cdot\left(\mathrm{tr}(\Ric^T)^2-|\Ric^T|^2\right)\\
&=&\frac{f^2}{|\nabla f|^2}\cdot\left(\frac{S^2}{4}-\frac{|\Ric(\eta)|^2}{|\eta|^2}-|\Ric_{\{\eta,\nu\}^\perp}|^2\right)\\
&=&\frac{f^2}{|\nabla f|^2}\cdot\left(\frac{S^2}{4}-\frac{S^2}{4}-|\Ric_{\{\eta,\nu\}^\perp}|^2\right)\\
&=&-\frac{f^2}{|\nabla f|^2}\cdot|\Ric_{\{\eta,\nu\}^\perp}|^2,
\end{eqnarray*}
from which $\Ric_{\{\eta,\nu\}^\perp}=0$ follows.

We have now all we need to conclude that both distributions $\mathrm{Span}(\eta,\nu)$ and its orthogonal complement are integrable and totally geodesic, the first one being of constant curvature $\epsilon$ and the second one being Ricci-flat (hence flat if $n=3$ or $4$).
Namely we already know that $\mathrm{Span}(\eta,\nu)$ is integrable since $[\eta,\nu]=0$.
Moreover,
\begin{eqnarray*}
\nabla_\eta\eta&=&\epsilon f\nabla f=\epsilon f|\nabla f|\nu\\
\nabla_\eta\nu&=&-\frac{f}{|\nabla f|}\Ric(\eta)=-\frac{Sf}{2|\nabla f|}\eta=-\frac{\epsilon f}{|\nabla f|}\eta\\
\nabla_\nu\eta&=&\nabla_\eta\nu=-\frac{\epsilon f}{|\nabla f|}\eta\\
\nabla_\nu\nu&=&0,
\end{eqnarray*}
so that all above expressions lie in $\mathrm{Span}(\eta,\nu)$, in particular $\mathrm{Span}(\eta,\nu)$ is totally geodesic.
As for $\mathrm{Span}(\eta,\nu)^\perp$, we compute, for all $X,Y\in\Gamma\left(\mathrm{Span}(\eta,\nu)^\perp\right)$,
\begin{eqnarray*}
\langle\nabla_XY,\eta\rangle&=&-\langle Y,\nabla_X\eta\rangle\\
&\bui{=}{\rm(\ref{eq:nablaeta})}&-\frac{f}{|\nabla f|}\cdot\left(\ric(\eta,X)\langle\nu,Y\rangle-\langle\nu,X\rangle\ric(\eta,Y)\right)\\
&=&0
\end{eqnarray*}
and, using $\mathrm{Span}(\eta,\nu)^\perp=\ker(\Ric)$,
\begin{eqnarray*}
\langle\nabla_XY,\nu\rangle&=&-\langle Y,\nabla_X\nu\rangle\\
&=&\frac{f}{|\nabla f|}\ric(X,Y)\\
&=&0.
\end{eqnarray*}
It follows that $\nabla_XY\in\Gamma\left(\mathrm{Span}(\eta,\nu)^\perp\right)$, therefore this distribution is integrable and totally geodesic.
To compute the curvature of both integral submanifolds, we notice that, from the above computations, $R_{\eta,\nu}\nu=\Ric(\eta)=\epsilon \eta$ and $R_{X,\nu}\nu=0=R_{X,\eta}\eta$ for all $X\in\ker(\Ric)$, so that 
$$\frac{\langle R_{\eta,\nu}\nu,\eta\rangle}{|\eta|^2}=\frac{\epsilon|\eta|^2}{|\eta|^2}=\epsilon$$
and, using the Gau\ss{} formula for curvature, for all $X,Y\in\Gamma\left(\mathrm{Span}(\eta,\nu)^\perp\right)$,
$$\ric_\Sigma(X,Y)=\ric(X,Y)-\langle R_{X,\frac{\eta}{|\eta|}}\frac{\eta}{|\eta|},Y\rangle-\langle R_{X,\nu}\nu,Y\rangle=0,$$
where we denoted by $\ric_\Sigma$ the Ricci curvature of the integral submanifold $\Sigma$ of $\mathrm{Span}(\eta,\nu)^\perp$.
Therefore, $\Sigma$ is Ricci-flat and thus flat if $1$-or $2$-dimensional.
% It remains to notice that, for all $X,Y\in\Gamma\left(\mathrm{Span}(\eta,\nu)^\perp\right)$,
% $$R_{X,Y}\eta=R_{X,Y}\nu=R_{X,\eta}\eta=R_{X,\nu}\nu=R_{\eta,X}Y=R_{\nu,X}Y=R_{\eta,\nu}X=0,$$
% that is, all ``mixed'' curvature terms vanish. USELESS SINCE WE ALREADY KNOW THAT $\mathrm{Span}(\eta,\nu)$ AND ITS ORTHOGONAL COMPLEMENT ARE INTEGRABLE AND TOTALLY GEODESIC
On the whole, this shows that the holonomy group of $M$ splits locally, therefore the universal cover of $M$ is isometric to the Riemannian product $S(\epsilon)\times\tilde{\Sigma}$ of the simply-connected complete surface with curvature $\epsilon\in\{-1,0,1\}$ with some simply-connected Ricci-flat manifold $\tilde{\Sigma}$.
In case $\epsilon=1$, the lift $\tilde{f}$ of $f$ to $\mathbb{S}^2\times\tilde{\Sigma}$ is constant along the $\tilde{\Sigma}$-factor and satisfies the equation $(\nabla^{\mathbb{S}^2})^2f=-f\cdot\mathrm{Id}$, which is exactly the equation characterizing the eigenfunctions associated to the first positive Laplace eigenvalue \cite[Theorem A]{Obata62}.
Furthermore, the isometry group of $\mathbb{S}^2\times\tilde{\Sigma}$ embeds into the product group of both isometry groups of $\mathbb{S}^2$ and $\tilde{\Sigma}$ and the first factor must be trivial since $\tilde{f}$, as the restriction of a linear form from $\R^3$ onto $\mathbb{S}^2$, is not invariant under $\{\pm\mathrm{Id}\}$.
Therefore, $M$ is isometric to $\mathbb{S}^2\times\Sigma$ for some Ricci-flat $\Sigma$ and $f$ is the trivial extension of an eigenfunction associated to the first positive Laplace eigenvalue on $\mathbb{S}^2$.\\
{ In case $\epsilon=0$, the manifold $M$ is Ricci-flat and therefore is isometric to the Riemannian product of $\R$ with a Ricci-flat manifold $\overline{N}$ as we mentioned in the introduction; our supplementary assumptions only mean that $\overline{N}$ carries a nontrivial parallel vector field. 
Mind in particular that $\overline{N}$ is not necessarily isometric to the Riemannian product of $\R$ or $\mathbb{S}^1$ with some Ricci-flat manifold, even if this is obviously locally the case.\\
In case $\epsilon=-1$, the lift $\tilde{f}$ of $f$ to $\mathbb{H}^2\times\tilde{\Sigma}$ is constant along the $\tilde{\Sigma}$-factor and satisfies the equation $(\nabla^{\mathbb{H}^2})^2f=f\cdot\mathrm{Id}$, which is exactly the Tashiro equation.
Since the isometry group of $\mathbb{H}^2\times\tilde{\Sigma}$ embeds into the product group of both isometry groups of $\mathbb{H}^2$ and $\tilde{\Sigma}$ and the first factor must be trivial since $\tilde{f}$ has no nontrivial symmetry \cite[Theorem 2 p.252]{Tashiro65}, we can deduce as above that $M$ is isometric to $\mathbb{H}^2\times\Sigma$ for some Ricci-flat $\Sigma$ and $f$ is the trivial extension of a solution to the Tashiro equation on $\mathbb{H}^2$.
This concludes the proof of Theorem \ref{t:eqfRicwitheta}.
}
\findemo 

% \begin{erems}\label{r:eqfRicwitheta}
% \noindent{\rm\begin{enumerate}
% \item 
% What happens if the assumption the restriction of $\eta$ to each level hypersurface of $f$ is parallel is dropped?
% \item Is the scalar curvature of $M$ always constant in case $M$ is closed?
% Up to now, no compact example with nonconstant scalar curvature could be constructed.
% \end{enumerate}
% }
% \end{erems}

{\bf Acknowledgment:}
%Part of this work was done while the first-named author enjoyed the friendly hospitality of the University of Greifswald and he would like to thank warmly Ines Kath for her generous support.
Part of this work was done while the second-named author received the support of the Humboldt Foundation which he would like to thank.
We also thank William Wylie for interesting discussions related to \cite{HePetersenWylie11102455} and \cite{HePetersenWylie11102456}.

\end{document}